\documentclass[final,leqno,onefignum,onetabnum]{siamltex1213}

\usepackage{graphicx,epsfig,color}
\usepackage{booktabs}
\usepackage[notref,notcite]{showkeys}
\usepackage{amsmath,cases}
\usepackage{amssymb,bm}%,mathdots}
\usepackage{multirow}%,pxfonts}
\newtheorem{remark}[theorem]{Remark}
\input{siam.sty}

\title{Any admissible harmonic Ritz value set is possible for prescribed GMRES residual norms}

\author{Kui Du\thanks{%Corresponding author. 
School of Mathematical Sciences and Fujian Provincial Key Laboratory of Mathematical Modelling and High-Performance Scientific Computation, Xiamen University, Xiamen 361005, China ({kuidu@xmu.edu.cn}). The research of this author was supported by the National Natural Science Foundation of China (No.11201392 and No.91430213), the Doctoral Fund of Ministry of Education of China (No.20120121120020), and the Fundamental Research Funds for the Central Universities (No.20720160002).} }

\begin{document}
\maketitle

\begin{abstract} %The roots of the polynomials GMRES generates to compute its residuals are Harmonic Ritz values. 
We show that any admissible harmonic Ritz value set is possible for prescribed GMRES residual norms, which is a complement for the results in [Duintjer Tebbens and Meurant, {\it SIAM J. Matrix Anal. Appl.}, 33 (2012), no. 3, pp. 958--978].  
\end{abstract}



\pagestyle{myheadings}
\thispagestyle{plain}
\markboth{Kui Du}{Admissible harmonic Ritz value set for GMRES}

\section{Introduction} The generalized minimal residual method (GMRES) \cite{saad1986gmres} is a popular iterative technique for solving large non-Hermitian linear systems. Greenbaum and Strako\u{s} \cite{greenbaum1994matri} proved that any convergence curve for the residual norm can be generated by GMRES applied to a non-derogatory matrix having prescribed eigenvalues. Greenbaum, Pt\'ak, and Strako\u{s} \cite{greenbaum1996any} showed later that any nonincreasing convergence curve is possible for GMRES. Arioli, Pt\'ak, and Strako\u{s} \cite{arioli1998krylo} gave a complete parametrization for the class of matrices having the same GMRES convergence curve. 

Recently, Duintjer Tebbens and Meurant \cite{tebbens2012any} showed that any Ritz value behavior is possible for prescribed GMRES residual norms. Since the roots of the polynomials GMRES generates to compute its residuals are harmonic Ritz values \cite{freund1992quasi}, it is interesting to investigate the harmonic Ritz value behavior for prescribed GMRES residual norms. 
Duintjer Tebbens and Meurant wrote \cite[page 974]{tebbens2012any} it is not clear whether any harmonic Ritz value behavior is possible for prescribed GMRES residual norms.
In this note we show that any admissible harmonic Ritz value set (introduced in \S 3) is possible for prescribed GMRES residual norms. 

The rest of this note is organized as follows. In the remainder of this section we introduce some notations. In section 2 we give some properties for GMRES and harmonic Ritz values. In section 3 we provide the main result of this note by exploiting a parameterized inverse eigenvalue problem.

To facilitate the discussion, we shall adopt the following notations. For a matrix ${\bf A}$, let $a_{ij}$, ${\bf a}_k$, $\tr({\bf A})$, $\det(\bf A)$, ${\bf A}_k$, $\wt{\bf A}_k$ and ${\bf A}^*$ denote the $i,j$ entry, the $k$th column, the trace, the determinant, the $k\times k$ principal submatrix, the matrix build with the first $k$ columns, and the conjugate transpose of $\bf A$,  respectively. The complex conjugate of a scalar $z$ is written $\wb z$. Let ${\bf e}_k$ denote the $k$th column of the identity matrix of appropriate order.

\section{Preliminaries} Let a nonsingular matrix ${\bf A}\in\mbbc^{n\times n}$ and a vector ${\bf b}\in\mbbc^n$ be given.  For an initial guess ${\bf x}_0$, GMRES approximates the exact solution of ${\bf Ax=b}$ at step $k$ by the vector ${\bf x}_k\in{\bf x}_0+\mcalk_k({\bf A},{\bf r}_0)$ that minimizes the Euclidean norm of the residual ${\bf r}_k:={\bf b}-{\bf Ax}_k$, i.e., $$\|{\bf r}_k\|=\min_{{\bf x}\in{\bf x}_0+\mcalk_k({\bf A},{\bf r}_0)}\|{\bf b-Ax}\|,$$ where %$\mbbp_k$ denotes the set of polynomials of degree $\leq k$ and 
$$\mcalk_k({\bf A},{\bf r}_0):=\spa\{{\bf r}_0,{\bf Ar}_0,\ldots,{\bf A}^{k-1}{\bf r}_0\}.$$
%Let $p_k(z)\in\mbbp_k$ with $p_k(0)=1$ such that $${\bf r}_k=p_k({\bf A}){\bf r}_0.$$
%%Let $${\bf K}_k=\bem\begin{array}{cccc}{\bf r}_0 & {\bf Ar}_0& \cdots & {\bf A}^{k-1}{\bf r}_0\end{array}\eem.$$ 
%\proposition
%The {\rm GMRES} polynomial $p_k(z)$ is given by $$p_k(z)=\prod_{j=1}^k\l(1-\frac{z}{\theta_j^{(k)}}\r).$$
%\endproposition

The Arnoldi process \cite{arnoldi1951princ} constructs an orthonormal basis of $\mcalk_k({\bf A},{\bf r}_0)$. Without loss of generality we assume $$\|{\bf r}_0\|=1,$$ and we also assume that the Arnoldi process for the pair $\{{\bf A},{\bf r}_0\}$ do not break down before the $n$th iteration. 
%This means the matrix $\bf A$ must be non-derogatory. 
Then we have the Arnoldi relation \beq\label{arnoldi}{\bf AV}={\bf V}{\bf H},\eeq where $\bf V$ is unitary and $\bf H$ is irreducible upper Hessenberg. The columns of $\wt{\bf V}_k$ form an orthonormal basis of $\mcalk_k({\bf A},{\bf r}_0)$.

%\subsection{Ritz values and harmonic Ritz values}
%The $k$ eigenvalues $\{\theta_j^{(k)}\}_{j=1}^k$ of ${\bf H}_k$ are called Ritz values {\rm(}at step $k${\rm)} of ${\bf A}$.  %The problem reduces to find a vector ${\bf y}_k\in\mbbc^k$ such that \beqas\|{\bf r}_0-{\bf AV}_k{\bf y}_k\|&=&\min_{{\bf y}\in\mbbc^k}\|{\bf r}_0-{\bf AV}_k{\bf y}\| = \min_{{\bf y}\in\mbbc^k}\l\|{\bf r}_0-{\bf V}_{k+1}{\bf\wt H}_{k}{\bf y}\r\|\\ &=& \min_{{\bf y}\in\mbbc^k}\l\|{\bf V}_{k+1}\l(\|{\bf r}_0\|{\bf e}_1-{\bf\wt H}_{k}{\bf y}\r)\r\|\\ &=& \min_{{\bf y}\in\mbbc^k}\l\|\|{\bf r}_0\|{\bf e}_1-{\bf\wt H}_{k}{\bf y}\r\|.\eeqas Therefore, $${\bf\wt H}_{k}^*{\bf\wt H}_{k}{\bf y}_k=\|{\bf r}_0\|{\bf\wt H}_{k}^*{\bf e}_1.$$ 
The eigenvalues of the generalized eigenvalue problem \beq\label{oGEP}{\bf \wt H}_k^*{\bf \wt H}_k{\bf z}=\theta{\bf H}_k^*{\bf z}\eeq are called harmonic Ritz values at step $k$ of the Arnoldi process for $\{{\bf A},{\bf r}_0\}$, giving the $k$-tuple $$\Theta^{(k)}=(\theta_1^{(k)},\theta_2^{(k)},\ldots,\theta_k^{(k)}).$$ For simplicity we assume these number are sorted in nondecreasing order (in magnitude). Note that $\theta_j^{(k)}$ is either nonzero finite complex number or $\infty$.  We denote by $\Theta$ the set \beq\label{Theta}\Theta:=\{\Theta^{(1)},\Theta^{(2)},\ldots,\Theta^{(n)}\}\eeq representing all $(n+1)n/2$ harmonic Ritz values.

%Note that the matrix ${\bf \wt H}_k$ we used is slightly different to that used in the literatures. 
%
%
%\proposition[Invariance under unitary similarity transformations]\label{invariance} Let  {\rm GMRES} be applied to the pair $\{{\bf A},{\bf r}_0\}$. If $\bf A$ is changed to $\bf U^*AU$ for some unitary matrix $\bf U$, and ${\bf r}_0$ is changed to ${\bf U^*r}_0$, the residual ${\bf r}_k$ change to ${\bf U^*}{\bf r}_k$. The Ritz values $\{\theta_j^{(k)}\}_{j=1}^k$ and the harmonic Ritz values $\{\theta_j^{(k)}\}_{j=1}^k$ do not change.
%\endproposition
%
%%See {\rm \cite{trefethen1997numer}} for a proof. 
%% Assume that {\rm GMRES} for $\{{\bf A},{\bf r}_0\}$ has not converged at step $n-1$, i.e., ${\bf r}_{n-1}\neq\bf 0$. 
% 
%% For any unitary matrix $\bf U$, the residual norms  of {\rm GMRES} for $\{{\bf U^*AU},{\bf U^*r}_0\}$ are the same as those of {\rm GMRES} for $\{{\bf H},\|{\bf r}_0\|{\bf e}_1\}$.   
%\remark Based on Proposition {\rm\ref{invariance}}, 
%the convergence curve generated by {\rm GMRES} applied to $\{{\bf A},{\bf r}_0\}$ is identical with the convergence curve generated by {\rm GMRES} applied to $\{{\bf H},{\bf e}_1\};$ the essentially all information about Ritz values, harmonic Ritz values, and residual norms is contained in $\bf H$.
%\endremark

Consider a QR factorization \beq\label{QRfactorization}{\bf H}={\bf Q}{\bf R},\eeq where $\bf Q$ is unitary irreducible upper Hessenberg and $\bf R$ is nonsingular upper triangular. By (\ref{QRfactorization}), we have the following factorizations of ${\bf H}_k$ and $\wt{\bf H}_k$: \beq\label{QRs}{\bf H}_k={\bf Q}_k{\bf R}_k,\qquad \wt{\bf H}_k=\wt{\bf Q}_k{\bf R}_k.\eeq Entries of $\bf Q$ and the relation to GMRES residual norms have been shown in \cite{meurant2012gmres,meurant2014neces}. For convenience of our investigation, we list some known results in Propositions \ref{matrixQ}, \ref{residual}, \ref{stagstructure}, and \ref{admissiblehrv}. We also give proofs for completeness.

\proposition\label{matrixQ} Rows $2$ through $n$ of the unitary irreducible upper Hessenberg matrix $\bf Q$ are uniquely  determined {\rm(}up to complex signs{\rm)} by the first row of $\bf Q$. Specifically, for $i=1:n-1$ and $j=i+1:n$
\beqas &&q_{i+1,i}=\rho_i\frac{\sqrt{1-\sum_{l=1}^i|q_{1l}|^2}}{\sqrt{1-\sum_{l=1}^{i-1}|q_{1l}|^2}}, \\ 
&&q_{i+1,j}=-\rho_i\frac{\wb q_{1i}q_{1j}}{\sqrt{1-\sum_{l=1}^{i-1}|q_{1l}|^2}\sqrt{1-\sum_{l=1}^i|q_{1l}|^2}},\eeqas where $$|\rho_1|=|\rho_2|=\cdots=|\rho_{n-1}|=1.$$
\endproposition

The proof of Proposition \ref{matrixQ} is straightforward by explicit calculations (exploiting the unitary irreducible upper Hessenberg structure of $\bf Q$). 

\proposition\label{residual} {\rm GMRES} residual norm at step $k$, $\|{\bf r}_k\|$, is given by  $$\|{\bf r}_k\|=\l(\sum_{l=k+1}^n |q_{1l}|^2\r)^{1/2}.$$
\endproposition
\proof It follows from $\wt{\bf H}_k=\wt{\bf Q}_k{\bf R}_k$ and \beqas\|{\bf r}_k\|&=&\min_{{\bf x}\in{\bf x}_0+\mcalk_k({\bf A},{\bf r}_0)}\|{\bf b-Ax}\| 
\\&=&\min_{{\bf y}\in\mbbc^k}\|{\bf r}_0-{\bf A}\wt{\bf V}_k{\bf y}\|
\\ &=& \min_{{\bf y}\in\mbbc^k}\|{\bf r}_0-{\bf V}{\bf \wt H}_k{\bf y}\| %\quad (\mbox {by } (\ref{arnoldi}))
\\ &=& \min_{{\bf y}\in\mbbc^k}\|{\bf e}_1-{\bf\wt H}_{k}{\bf y}\| %\quad (\mbox{by } {\bf V} unitary and r_0=Ve_1)
%\\&=&\|{\bf e}_1-\wt{\bf H}_k(\wt{\bf H}_k^*\wt{\bf H}_k)^{-1}\wt{\bf H}_k^*{\bf e}_1\|,
\eeqas
%It follows from $$\|{\bf r}_k\|=\min_{{\bf x}\in{\bf x}_0+\mcalk_k({\bf A},{\bf r}_0)}\|{\bf b-Ax}\| =\min_{{\bf y}\in\mbbc^k}\|{\bf r}_0-{\bf A}\wt{\bf V}_k{\bf y}\|,$$ the Arnoldi relation ${\bf A}\wt{\bf V}_k={\bf V}\wt{\bf H}_k$, and the QR factorization $\wt{\bf H}_k=\wt{\bf Q}_k{\bf R}_k$ that
%\beqas\min_{{\bf x}\in{\bf x}_0+\mcalk_k({\bf A},{\bf r}_0)}\|{\bf b-Ax}\| &=&\min_{{\bf y}\in\mbbc^k}\|{\bf r}_0-{\bf A}\wt{\bf V}_k{\bf y}\|= \min_{{\bf y}\in\mbbc^k}\|{\bf r}_0-{\bf V}{\bf \wt H}_k{\bf y}\|\\ &=& \min_{{\bf y}\in\mbbc^k}\|{\bf V}({\bf e}_1-{\bf\wt H}_{k}{\bf y})\|=\min_{{\bf y}\in\mbbc^k}\|{\bf e}_1-{\bf\wt H}_{k}{\bf y}\|\\&=&\|{\bf e}_1-\wt{\bf H}_k(\wt{\bf H}_k^*\wt{\bf H}_k)^{-1}\wt{\bf H}_k^*{\bf e}_1\|,\eeqas 
that
\beqas\|{\bf r}_k\|=\|({\bf I}-{\bf\wt Q}_{k}{\bf \wt Q}_{k}^*){\bf e}_1\|=\l(\sum_{l=k+1}^n|{\bf q}_{l}^*{\bf e}_1|^2\r)^{1/2}=\l(\sum_{l=k+1}^n |q_{1l}|^2\r)^{1/2}. \quad \endproof\eeqas  

Proposition \ref{residual} implies that the {\rm GMRES} residual norms can be read from the first row of $\bf Q$. 
%The following proposition is a direct result of Propositions \ref{matrixQ} and \ref{residual}.

\proposition\label{stagstructure} {\rm GMRES} applied to  $\{{\bf A},{\bf r}_0\}$ stagnates at step $k$, i.e., $$\|{\bf r}_k\|=\|{\bf r}_{k-1}\|,$$ if and only if ${\bf q}_k$, the $k$th column of $\bf Q$, satisfies $${\bf q}_k=q_{k+1,k}{\bf e}_{k+1}.$$
\endproposition
\proof It follows from Proposition \ref{residual} and $\|{\bf r}_k\|=\|{\bf r}_{k-1}\|$ that  $q_{1k}=0.$ Then $q_{ik}=0$ for $i=2:k$ follows from Proposition \ref{matrixQ}.  Therefore ${\bf q}_k=q_{k+1,k}{\bf e}_{k+1}$. Conversely, if ${\bf q}_k=q_{k+1,k}{\bf e}_{k+1}$, by Proposition \ref{residual},  $\|{\bf r}_k\|=\|{\bf r}_{k-1}\|$.
\endproof

Next, we characterize the harmonic Ritz values when GMRES stagnates.  By (\ref{QRs}),  the generalized eigenvalue problem (\ref{oGEP}) for harmonic Ritz values at step $k$ reduces to: \beq\label{gep}{\bf R}_k{\bf z}=\theta{\bf Q}_k^*{\bf z}.\eeq

\proposition\label{admissiblehrv} Assume that {\rm GMRES} applied to $\{{\bf A},{\bf r}_0\}$ stagnates from step $k+1$ to step $k+m$ {\rm(}$0\leq k< k+m \leq n-1${\rm)}, i.e., $$\|{\bf r}_k\|=\|{\bf r}_{k+1}\|=\cdots=\|{\bf r}_{k+m}\|.$$ Then harmonic Ritz values $\{\theta_j^{(k+i)}\}_{j=1}^{k+i}$ at step $k+i$ for $i=1:m$ satisfy $$\theta_j^{(k+i)}=\l\{\begin{array}{ll} \theta_j^{(k)}, & 1\leq j\leq k \\ \infty, & k+1\leq j\leq k+i.\end{array}\r.$$ 
\endproposition

\proof By Proposition \ref{stagstructure}, the generalized eigenvalue problem  (\ref{gep}) at step $k+i$ ($1\leq i\leq m$) reduces to   
$$ {\bf R}_{k+i}{\bf z}= \theta \bem {\bf Q}_k^* & \wb q_{k+1,k}{\bf e}_k{\bf e}_1^*\\ {\bf 0} & {\bf T} \eem \bf z. $$ The statement follows from ${\bf R}_{k+i}$ is nonsingular upper triangular, and all diagonal entries of the upper triangular matrix ${\bf T}$ are zero.
\endproof

\section{Harmonic Ritz values for prescribed GMRES residual norms} We call a set $\Theta$ defined in $(\ref{Theta})$ satisfying $\infty\notin\Theta$ an admissible harmonic Ritz value  set for stagnation-free GMRES. We also call a set $\Theta$ defined in $(\ref{Theta})$ satisfying Proposition ${\ref{admissiblehrv}}$ an admissible harmonic Ritz value set for {\rm GMRES} with stagnation.  In this section we will show that any admissible harmonic Ritz value set is possible for given GMRES residual norms.

By Propositions \ref{matrixQ} and \ref{residual}, given GMRES residual norms implies entries of $\bf Q$ are uniquely determined up to complex signs. Given an admissible harmonic Ritz value set $\Theta$, by constructing a desirable nonsingular upper triangular matrix $\bf R$ (see the approach below), we can obtain a pair $\{{\bf H},{\bf e}_1\}$, for which GMRES produces harmonic Ritz value set $\Theta$ and the prescribed residual norms.  

Now we describe how to construct $\bf R$. At step $k$, let $\{\wh r_{jk}\}_{j=1}^k$ denote the entries of the last column of ${\bf R}_k^{-1}$. 

(i) If ${\bf Q}_k$ is nonsingular, we consider the parameterized inverse eigenvalue problem for prescribed harmonic Ritz values: Given $k$ nonzero complex number (harmonic Ritz values) $\theta_1^{(k)},\theta_2^{(k)},\ldots,\theta_k^{(k)}$, and two matrices ${\bf Q}_k$, ${\bf R}_{k-1}^{-1}$;  find $\{\wh r_{jk}\}_{j=1}^k$ such that $$\l\{\frac{1}{\theta_1^{(k)}}, \frac{1}{\theta_2^{(k)}}, \ldots, \frac{1}{\theta_k^{(k)}}\r\}$$ is the spectrum of the matrix $${\bf Q}_k^*{\bf R}_k^{-1}={\bf Q}_k^*\bem{\bf R}_{k-1}^{-1} & {\bf 0}\\ {\bf 0} & 0\eem+\sum_{j=1}^k\wh r_{jk}{\bf Q}_k^*{\bf e}_j{\bf e}_k^*.$$ Note that the $k$ matrices ${\bf Q}_k^*{\bf e}_j{\bf e}_k^*$ $(1\leq j\leq k)$ are linearly independent and $$\tr({\bf Q}_k^*{\bf e}_1{\bf e}_k^*)=q_{1,k}\neq 0.$$ Helton et. al \cite{helton1997matri} proved that almost all such parameterized inverse eigenvalue problems are solvable. 

(ii) If ${\bf Q}_k$ is singular, which means GMRES stagnates at step $k$, by Proposition \ref{admissiblehrv}, we can set $\wh r_{jk}=0$ for $j=1:k-1$ and $\wh r_{kk}=1$. 

Once all $\{\wh r_{jk}\}_{j=1}^k$ for $k=1:n$ are found,  we obtain ${\bf R}^{-1}$, then $\bf H$ follows by ${\bf H=QR}$. GMRES applied to $\{{\bf H},{\bf e}_1\}$ produces the prescribed harmonic Ritz values (in all steps) and the prescribed residual norms. That is to say, any admissible harmonic Ritz value set is possible for given GMRES residual norms.

\remark Given $\bf Q$ and ${\bf R}_{k-1}$, we provide an obvious approach for construction of ${\bf R}_k$. It is sufficient to consider the case ${\bf Q}_k$ is nonsingular.  We obtain $\{r_{jk}\}_{j=1}^k$, the entries of the last column of ${\bf R}_k$, by solving the following linear system  $$\det({\bf R}_k-\theta_i^{(k)}{\bf Q}_k^*)=0,\quad i=1:k.$$ 

\endremark

%\cite{tebbens2014presc}
%\cite{meurant2015role}
%\cite{greenbaum1996any}
%\cite{arioli1998krylo}
%\cite{cao1997note}

%\bibliographystyle{siam}
%\bibliography{/users/dukui/mywork/common/macbib}

\begin{thebibliography}{10}

\bibitem{arioli1998krylo}
{\sc M.~Arioli, V.~Pt{{\'a}}k, and Z.~Strako{\v{s}}}, {\em Krylov sequences of
  maximal length and convergence of {GMRES}}, BIT, 38 (1998), pp.~636--643.

\bibitem{arnoldi1951princ}
{\sc W.~E. Arnoldi}, {\em The principle of minimized iteration in the solution
  of the matrix eigenvalue problem}, Quart. Appl. Math., 9 (1951), pp.~17--29.

\bibitem{tebbens2012any}
{\sc Jurjen Duintjer~Tebbens and G{{\'e}}rard Meurant}, {\em Any {R}itz value
  behavior is possible for {A}rnoldi and for {GMRES}}, SIAM J. Matrix Anal.
  Appl., 33 (2012), pp.~958--978.

\bibitem{freund1992quasi}
{\sc Roland~W. Freund}, {\em Quasi-kernel polynomials and their use in
  non-{H}ermitian matrix iterations}, J. Comput. Appl. Math., 43 (1992),
  pp.~135--158.

\bibitem{greenbaum1996any}
{\sc Anne Greenbaum, Vlastimil Pt{{\'a}}k, and Zden{\v{e}}k Strako{\v{s}}},
  {\em Any nonincreasing convergence curve is possible for {GMRES}}, SIAM J.
  Matrix Anal. Appl., 17 (1996), pp.~465--469.

\bibitem{greenbaum1994matri}
{\sc Anne Greenbaum and Zden{\v{e}}k Strako{\v{s}}}, {\em Matrices that
  generate the same {K}rylov residual spaces}, in Recent advances in iterative
  methods, vol.~60 of IMA Vol. Math. Appl., Springer, New York, 1994,
  pp.~95--118.

\bibitem{helton1997matri}
{\sc William Helton, Joachim Rosenthal, and Xiaochang Wang}, {\em Matrix
  extensions and eigenvalue completions, the generic case}, Trans. Amer. Math.
  Soc., 349 (1997), pp.~3401--3408.

\bibitem{meurant2012gmres}
{\sc G{{\'e}}rard Meurant}, {\em G{MRES} and the {A}rioli, {P}t{\'a}k, and
  {S}trako\v s parametrization}, BIT, 52 (2012), pp.~687--702.

\bibitem{meurant2014neces}
\leavevmode\vrule height 2pt depth -1.6pt width 23pt, {\em Necessary and
  sufficient conditions for {GMRES} complete and partial stagnation}, Appl.
  Numer. Math., 75 (2014), pp.~100--107.

\bibitem{saad1986gmres}
{\sc Youcef Saad and Martin~H. Schultz}, {\em G{MRES}: a generalized minimal
  residual algorithm for solving nonsymmetric linear systems}, SIAM J. Sci.
  Statist. Comput., 7 (1986), pp.~856--869.

\end{thebibliography}

\end{document}